\documentclass[twoside]{article}

\usepackage[preprint]{aistats2027}
\usepackage{amsmath,amsfonts,bm}

\def\1{\bm{1}}

\def\mK{{\bm{K}}}

\def\mQ{{\bm{Q}}}

\def\mV{{\bm{V}}}

\DeclareMathAlphabet{\mathsfit}{\encodingdefault}{\sfdefault}{m}{sl}
\SetMathAlphabet{\mathsfit}{bold}{\encodingdefault}{\sfdefault}{bx}{n}

\usepackage{algorithm}
\usepackage{algpseudocode}
\algrenewcommand\algorithmicdo{}
\usepackage{hyperref}
\usepackage{url}
\usepackage{graphicx,epstopdf}
\usepackage{dcolumn}
\usepackage{bm}
\usepackage{amsmath,amssymb,tabularx,booktabs}
\usepackage{amsthm}
\usepackage{cuted}
\usepackage{enumerate}
\usepackage{stfloats}
\usepackage{enumitem}
\usepackage{subcaption}
\usepackage[dvipsnames]{xcolor}
\usepackage{float}
\usepackage[skip=0.5ex]{subcaption}
\usepackage{lipsum}
\usepackage{framed}
\usepackage[most]{tcolorbox}
\usepackage{balance}

\newtheorem*{lemma*}{Lemma}
\newtheorem*{theorem*}{Theorem}
\newtheorem{theorem}{Theorem}

\theoremstyle{definition}

\theoremstyle{remark}
\newtheorem{remark}[theorem]{Remark}

\usepackage[round]{natbib}

\hypersetup{
    colorlinks=true,
    linkcolor=blue,
    filecolor=magenta,      
    urlcolor=blue,
    citecolor=blue,
    pdftitle={Overleaf Example},
    pdfpagemode=FullScreen,
    }
\begin{document}

%

%

\twocolumn[

\aistatstitle{Mean-field Token Steering}

\aistatsauthor{ Md Zulfiqur Haider \And Asmaa Eldesoukey \And  Abhishek Halder \And Yongxin Chen }

\aistatsaddress{ Iowa State University \And  Iowa State University \And Iowa State University \And Georgia Institute of Technology } ]

\begin{abstract}
We consider stochastic token steering through a deep stack of transformer layers as a mean-field control problem for an interacting particle system on the unit sphere, where interactions arise from the attention mechanism. In the mean-field limit, the token population evolves via a controlled integro-partial differential equation where the feedback control acts as a corrective vector field guiding the population toward an aggregate configuration of tokens in the representation space. We deduce the minimum effort control to accomplish this task 
over a fixed time horizon interpreted as normalized layer depth. We show that the resulting controlled dynamics for the token population can be understood as a mean-field Schrödinger bridge adapted to the sphere geometry. Building on this, we propose a novel generalization of the dynamic Sinkhorn recursion to compute the consistent token evolution. Numerical experiments confirm that the proposed approach achieves exact steering of the token population.
\end{abstract}

\section{Introduction}

The transformer architecture \citep{vaswani2017attention} and the self-attention mechanism \citep{bahdanau2014neural,luong2015effective} within it are now widely used in large language models (LLMs). An important reason behind this wide usage and the significant impact made thereby lies in their ability to capture long-range dependencies in data, such as in prompts with a large number of tokens. Notwithstanding the empirical success, the mathematical understanding of how representations evolve through a deep stack of transformer layers via interactions among the tokens is only beginning to emerge now. Still, there remains a significant theoretical gap in understanding what kind of representations can emerge through such evolutions.

One promising approach toward this understanding is to index the layers by time $t$ and to consider the evolution of $N$ interacting tokens in the large-depth limit as a continuous-time dynamical system. Then, the token ensemble evolves as a system of $N$ interacting particles, where a particle represents an individual token embedding \citep{lu2019transformerode,dutta2021redesigning,Sander2022,Geshkovskietal025}. Mathematically, this evolution is described by a system of $N$ coupled ordinary or stochastic differential equations, depending on whether the evolution is deterministic or noisy. Of particular interest is the regime of a large number of tokens ($N\rightarrow\infty$), referred to as the \emph{mean-field limit}, in which case the dynamics are governed by a nonlinear integro-partial differential equation (PDE) that describes the evolution of the \emph{token distribution over the unit sphere $\mathbb{S}^{d-1}$}, where $d$ is the embedding dimension. The spherical geometry arises from the layer normalization. This integro-PDE is parameterized by the query, key, and value matrices $\bm{Q}_t,\bm{K}_t,\bm{V}_t$ indexed by layer or time index $t$, and the inverse temperature $\beta>0$ in the self-attention mechanism. 
Understanding the mean-field limit is of the essence in practice as the context lengths and depths in LLMs continue to increase. For example, both Claude Opus 5.5 and Gemini frontier models feature 1-million-token context windows, while GPT-4o features a 128k-token context window.

Previous studies on the aforementioned integro-PDE found that the token dynamics can lead to the emergence of a single cluster in the long run \citep{geshkovski2023emergence,Geshkovskietal025}, and dynamic meta-stability with partial clustering at intermediate times \citep{geshkovski2024metastability,Bruno2025b,bruno2026multiscale}. In practice, the clustering relates to rank or entropy collapse \citep{dong2021attention,zhai2023stabilizing,giorlandino2026two}, which has been verified empirically as well as via non-mean-field theoretical studies \citep{markdahl2017almost,abella2025consensus,criscitiello2026synchronization}. 

While existing works on mean-field token dynamics have focused on \emph{analysis}, a natural question is whether it is possible to \emph{design} an optimal correction of the nonlocal self-attention dynamics to steer the token population toward a desired aggregate configuration. This work shows why and how this is possible. 

The optimal correction we derive can be interpreted as \emph{additive control} of the noisy mean-field interaction between tokens. By standard Lie-Trotter splitting \citep[p. 47]{hairer2006geometric}, the corrective drift associated with additive control can be realized through a feedforward layer concatenated to the self-attention layer \citep[Remark 2.1]{bruno2026multiscale}, \citep{isobe2026training}. To reduce computational cost, it is preferable to achieve such active steering in a \emph{fixed time horizon} $[0,T]$ interpreted as normalized layer depth. In this work, we put forth this new perspective of active guidance of the token distribution by shifting focus from mean-field dynamics to mean-field control.  

\paragraph{Contributions} We formulate and solve the mean-field token steering problem through the lens of stochastic optimal control. Our specific contributions are the following.
\begin{enumerate}
    \item We show that the minimum effort control that steers the interacting token ensemble from a given aggregate configuration to another can be understood as a mean-field Schrödinger bridge adapted to the sphere geometry. This establishes a new connection between two lines of literature: mean-field self-attention dynamics and Schr\"{o}dinger bridges. While Schr\"{o}dinger bridges and their Riemannian versions are well-known in diffusion-based generative AI \citep{de2021diffusion,shi2023diffusion,albergo2025stochastic}, the variant of the Schrödinger bridge appearing here is far from standard since it involves both nonlocal self-attention as well as the spherical geometry.

    \item Building on the above connection, we derive the mathematical structure of the minimum-effort control and the optimal evolution of the token distribution, in terms of the solution to a system of coupled nonlinear integro-PDE boundary value problems. We emphasize how these results significantly extend the previously known results for Schr\"{o}dinger bridges.

    \item To numerically solve the derived system of equations for the minimum-effort active token steering, we propose a novel generalization of the dynamic Sinkhorn recursion (also known as the iterative proportional fitting algorithm). We illustrate the proposed algorithm for steering tokens over $\mathbb{S}^{1}$ and $\mathbb{S}^{2}$.  
\end{enumerate}

The assumptions needed in our development are mild. For instance, we work with normalized self-attention instead of its unnormalized simplification \citep{rigollet2026mean,isobe2026training}. Also, we only require the matrices $\bm{Q}_t,\bm{K}_t,\bm{V}_t$ to be continuous in $t\in[0,T]$. In comparison, the mean-field transformer dynamics literature makes stronger assumptions on these matrices with varying degrees of generality. These include e.g., uniformly continuous and bounded $\bm{Q}_t^{\top}\bm{K}_t$ and $\bm{V}_t\equiv\bm{I}$ (identity matrix) in \citet{abella2025consensus}, time (i.e., depth) invariance as in $\bm{Q}_t\equiv\bm{Q},\bm{K}_t\equiv\bm{K},\bm{V}_t\equiv\bm{V}$ in \citet{geshkovski2023emergence}, time-invariant invertible $\bm{Q},\bm{K},\bm{V}$ in \citet{bruno2026multiscale}, $\bm{Q}^{\top}\bm{K}=\bm{V}$ symmetric in \citet{Sander2022,Bruno2025b,Geshkovskietal025,Burgeretal2025}; the latter is helpful for a Wasserstein-like gradient flow interpretation but is not needed in our setting.  

\paragraph{Related works} The idea of identifying layers with time and then deriving a continuous-time flow in the large-depth limit is reminiscent of neural ordinary and stochastic differential equations \citep{chen2018neural,tzen2019neural}. The transformer case differs from that non-interacting setting because long-range interactions among tokens give rise to a mean-field integro-PDE for the token distribution. The transformer then acts as a flow map on $\mathcal{P}(\mathbb{S}^{d-1})$, the space of probability distributions on $\mathbb{S}^{d-1}$. 

We note that a mean-field integro-PDE also arises \citep{mei2018mean,chizat2018global} in training a feedforward network with a single hidden layer of infinite width via stochastic gradient descent. Unlike that case, the transformer mean-field dynamics is not induced by training, but instead by the propagation of the distributional state of the tokens.    

In this work, we consider \emph{noisy} nonlocal interactions among tokens: the presence of noise makes the mean-field integro-PDE second-order \citep{agazzi2026stochastic}, i.e., a McKean-Vlasov-Fokker-Planck PDE \citep{mckean1966class,sznitman2006topics}. Prior literature on stochastic transformers has motivated the noise due to random initialization of the value matrix \citep{fedorov2026clustering}, or as computational regularization \citep{ICLR2024_4bdeeaeb,isobe2026training}.

\emph{Control in transformers} has been used from different viewpoints: as Low Rank Adaptation (LoRA) for the purpose of fine-tuning of pre-trained models \citep{zhang2024parameter}, as positional encoding and prompt prefixes to counteract mode collapse \citep{imaizumi2026anti}, and as training over a parametrized family of transformers \citep{kan2026optimal,isobe2026training}. Different from these works, we formulate mean-field token steering as a minimum-effort distributional control problem. In doing so, we make connections with the broader \emph{mean-field Schr\"{o}dinger bridge} literature \citep{backhoff2020mean,zhang2026modeling,eldesoukey2026generalized,eldesoukey2026schrodingerbridgeskineticswarming}.

\paragraph{Organization} In Section \ref{sec:formulation}, we introduce both the uncontrolled and controlled versions of the transformer mean-field dynamics. In Section \ref{sec:bridges}, we present our main results: the problem formulation, mathematical structure of the Schr\"{o}dinger system that determines the optimal controller for mean-field token steering, and a computational method to solve the same. Numerical experiments are reported in Section \ref{sec:Numerical}. Section \ref{sec:Conclusions} closes with concluding remarks. 


\section{transformer Mean-field Dynamics}
\label{sec:formulation}
To account for layer normalization, we consider token embeddings that evolve on $\mathbb{S}^{d-1}$. Following \citet{rigollet2026mean,shalova2026solutions,isobe2026training}, we model their evolution over a continuum of attention layer depth $t\in [0,T]$ as the It\^{o} stochastic differential equation:
\begin{align}
    {\rm d} X_t^{i} &= P_{X_t^i}^{\perp}\left( \frac{1}{\mathcal Z_{t}^i}\sum_{j=1}^N \Gamma_t(X_t^i, X_t^j)\, \mV_t \, X_t^j\right) {\rm d}t + \sqrt{\sigma}\, {\rm d}B_t^i, \label{eq:finite-interacting-transformer} 
\end{align}
where $X_t^i \in \mathbb{S}^{d-1}$ denotes the embedding (state) of the $i$-th token for $i \in [N]$, and $\mV_t \in \mathbb{R}^{d\times d}$ is the value matrix. 

The kernel $\Gamma_t : \mathbb{S}^{d-1} \times \mathbb{S}^{d-1} \to \mathbb{R}_+$ denotes the unnormalized attention given by the $i$-th token to the $j$-th token. The corresponding normalization $\mathcal Z_t^{i}:=\sum\limits_{j=1}^N \Gamma_t(X_t^i, X_t^j)$. Specifically, the \emph{unnormalized attention kernel} 
 \begin{align}
     \Gamma_t(x,y) = e^{\beta\langle \mQ_t x, \mK_t y\rangle_{\mathbb R^d}} \qquad \text{for}\;x,y\in\mathbb{S}^{d-1},
     \label{defGammat}
 \end{align}
where $\mQ_t,\mK_t$ are the query and key matrices at layer $t$, respectively, and $\beta \in \mathbb R_{++}$ is a 
scaling parameter that represents the inverse temperature. Both $\mQ_t,\mK_t$ are of size: sequence length$\times d$, so $\mQ_t^{\top}\mK_t \in \mathbb{R}^{d\times d}$. In \eqref{defGammat} and throughout, $\langle \cdot, \cdot \rangle_{\mathbb R^d}$ denotes the Euclidean inner product.

In \eqref{eq:finite-interacting-transformer}, the attention aggregated by the $i$-th token amounts to a drift in its own dynamics. The operator $P_{X_t^i}^{\perp}$ projects this drift onto the tangent space $T_{X_t^i} \mathbb S^{d-1}$, i.e., $$P_x^{\perp} (v)  :=  v - \langle x,v\rangle_{\mathbb{R}^{d}} \, x$$ for $v \in \mathbb R^d$. 
Furthermore, $B_t^i$ denotes a Brownian motion on $\mathbb{S}^{d-1}$, and $\sqrt{\sigma}$ denotes noise intensity. 

\paragraph{Assumptions} We assume that
\begin{itemize}
    \item[\textbf{A1.}] the matricial trajectories $t \mapsto \bm{Q}_t, t \mapsto \bm{K}_t, t \mapsto \bm{V}_t$ remain continuous over $t\in[0,T]$,
    \item[\textbf{A2.}] the Brownian motions $\{B_t^i\}_{i=1}^{N}$ are independent. 
\end{itemize}
\paragraph{Mean-field limit} To describe the evolution of the token embeddings as a \emph{population}, we consider the empirical measure 
\begin{align*}
    \rho_t^N := \frac{1}{N} \sum_{i=1}^{N} \delta_{X_t^i},
\end{align*}
where $\delta_{X_t^i}$ denotes the Dirac mass concentrated  at $X_t^i$.  In the mean-field limit $N \to \infty$, interactions among the tokens through the attention mechanism take place via averaged effects \citep{bolley2011stochastic}. When the initial token states are sampled independently from a distribution $\rho_{\rm initial}\in\mathcal{P}(\mathbb{S}^{d-1})$, 
then by standard propagation-of-chaos arguments \citep{sznitman2006topics,ChaintronandDiez022}, the empirical (random) measure $\rho_t^N$ converges to a deterministic probability distribution  $\rho_t\in\mathcal{P}(\mathbb{S}^{d-1})$ in almost sure sense.  Throughout, we use the same symbol $\rho_t$ to denote a probability law in $\mathcal P(\mathbb S^{d-1})$ and its probability density function
, i.e., $\rho_t ({\rm d}x) = \rho_t(x)\:{\rm d} x$.

As a result of propagation of chaos, \eqref{eq:finite-interacting-transformer} converges to the following McKean--Vlasov-type dynamics on $\mathbb{S}^{d-1}$ \citep{PASCUCCI2026104978}:
\begin{equation}
    {\rm d} X_t^i = P_{X_t^i}^{\perp}\Big( b_t[\rho_t](X_t^i)\Big)\, {\rm d}t + \sqrt{\sigma}\, {\rm d}B_t^i, \quad t\in[0,T],
\label{eq:mean-field-transformer}
\end{equation}
where 
\begin{align}
    b_t[\rho_t](x) &:= \frac{1}{{\mathcal Z_t[\rho_t](x)}} \int_{\mathbb S^{d-1}} \Gamma_t(x,y)\,\mV_t \, y\,\rho_t(y) \, {\rm d}y, \label{eq:mean-field-force} 
\end{align}
and the attention normalization
\begin{align}
    \mathcal Z_t[\rho_t](x) &:= \int_{\mathbb S^{d-1}}  e^{\beta\langle \mQ_t x,\mK_t y\rangle_{\mathbb R^d}}  \rho_t(y) \, {\rm d}y.
\end{align}
We view the nonlocal vector field $b_t$ as the attention-induced mean-field velocity for the token embeddings. In the control context, we will treat this as a prior dynamics resulting from self-attention. We proceed by dropping the index $i$ for brevity.

\paragraph{Controlled dynamics and mean-field integro-PDE} Recent work by \citet{isobe2026training} studies the competition between the token clustering effected by the attention kernel and the feed-forward term, viewing the latter as a \emph{trainable controller} that can steer tokens away from the clustered distribution. The present work builds on this perspective. Specifically, we consider the question of \emph{how to design a suitable controller that steers potentially clustering token dynamics from a prescribed initial distribution toward a desired target distribution}. Our analysis below explicitly identifies the optimal controller to accomplish the steering task.

To this end, we consider the controlled dynamics 
\begin{equation}\label{eq:dynamics-sphere-meanfield-controlled}
    {\rm d} X_t = P_{X_t}^{\perp}\Big(b_t[\rho_t](X_t) + u_t(X_t)\Big)\, {\rm d}t + \sqrt{\sigma}\, {\rm d}B_t,
\end{equation}
where $u_t: \mathbb S^{d-1} \to \mathbb R^{d}$ represents the steering control input.
Equivalently, the token distribution satisfies the mean-field integro-PDE
\begin{align}
    \partial_t \rho_t(x) + \nabla \cdot \Big(\rho_t (x)\, P_x^{\perp}\Big(b_t[\rho_t](x)& + u_t(x)\Big)\Big) \nonumber \\
    &=  \frac{\sigma}{2}\, \Delta \rho_t(x),
\label{ControlledMeanFieldPDE}    
\end{align}
where $\nabla \cdot (\cdot)$ denotes the spherical divergence, and $\Delta$ is the Laplace-Beltrami operator on $\mathbb{S}^{d-1}$. We interpret $u_t(X_t)$ as a corrective velocity field for the mean-field token dynamics.

\section{Schr\"odinger Bridges over Attention Dynamics: Formulation and Results}\label{sec:bridges}
\subsection{Problem Formulation}
We propose to cast the controlled mean-field token steering as a Schr\"odinger bridge problem, where the control design is informed by two probability laws $\rho_{\rm initial}, \rho_{\rm final} \in \mathcal P (\mathbb S^{d-1})$. We think of $\rho_{\rm initial}, \rho_{\rm final}$ as given initial and target token distributions.

Our goal is to identify a controlled evolution $\rho_t$ over $t\in [0,T]$ that respects the mean-field token distribution dynamics \eqref{ControlledMeanFieldPDE}, satisfies $\rho_0=\rho_{\rm initial}, \rho_T = \rho_{\rm final}$, and minimizes the total control effort 
\begin{align}
   \int_0^T  \int_{\mathbb{S}^{d-1}}  \|u_t(x)\|_{\mathbb R^d}^2 \, \rho_t(x)\, \mathrm{d}x\,\mathrm{d}t,
\label{TotalControlEffort}   
\end{align}
where $\Vert \cdot \Vert_{\mathbb R^d}^2$ denotes the Euclidean norm. We can interpret the spatial integral $$\int_{\mathbb{S}^{d-1}}  \|u_t(x)\|_{\mathbb R^d}^2 \, \rho_t(x)\, \mathrm{d}x$$ as the instantaneous kinetic energy for the token population, and the space-time integral \eqref{TotalControlEffort} as the total work done due to the population transfer induced by the control $u_t$ over $t\in[0,T]$.

Notice that while $P_x^\perp (u_t(x))$ appears in the token dynamics \eqref{eq:dynamics-sphere-meanfield-controlled}, any component of $u_t(x)$ normal to the tangent space of the sphere increases control energy cost. Thus, we can restrict ourselves to controllers that belong to the tangent space. We then let the set of admissible control trajectories be defined as
\begin{align}
\!\!\!    \mathcal  U (\rho_\cdot) &:= \Big \{[0,T] \ni t \mapsto u_t: u_t(x) \in T_x \mathbb S^{d-1} \text{ for } x \in \mathbb S^{d-1}, \nonumber \\
  &   \text{and} \int_0^T  \int_{\mathbb{S}^{d-1}}  \|u_t(x)\|^2  \, \rho_t(x)\, \mathrm{d}x\,\mathrm{d}t < \infty \Big\},
\end{align}
where $\Vert \cdot \Vert$ denotes the round metric norm.

Letting
\begin{align}
    \mathcal D: = \left \{[0,T] \ni t \mapsto \rho_t : \rho_t \in \mathcal P(\mathbb S^{d-1}) \right \},
\end{align}
we now state our main problem as follows.

\begin{subequations}
Determine \label{problem:spherical_control}
\begin{align}
& (u_\cdot^\star, \rho_\cdot^\star)    \in  \nonumber \\
 &\underset{(u_\cdot, \rho_\cdot) \,  \in  \, \mathcal U (\rho_\cdot)  \times \mathcal D}{\arg\:\min}
 \frac{1}{2 \sigma}     \int_0^T
    \int_{\mathbb{S}^{d-1}}   \!\! \|u_t(x)\|^2  \,     \rho_t(x)\,     \mathrm{d}x\,\mathrm{d}t
    \label{eq:spherical_cost}    
    \end{align}
    subject to
    \begin{align}
    &  \partial_t \rho_t(x) + \nabla \cdot \Big(\rho_t(x)\, \Big(P_x^{\perp} \big(b_t[\rho_t](x)\big) + u_t(x)\Big) \Big) \nonumber \\
    & \qquad \qquad \qquad \qquad \qquad \qquad= \frac{\sigma}{2}\, \Delta \rho_t(x),    \label{eq:spherical_FP}    \\
    &      \rho_0= \rho_{\rm initial},     \qquad     \rho_T=\rho_{\rm final}.    \label{eq:spherical_endpoints}
 \end{align}
\end{subequations}


In what follows, we present the structure of the solution for the constrained variational problem \eqref{problem:spherical_control}.

\subsection{Identification of the Schr\"odinger Bridge}
The minimizer for problem~\eqref{problem:spherical_control} comprises of the optimal control $u_t^{\star}$ and the optimally controlled $\rho_t^{\star}$ for $t\in[0,T]$. The optimal control policy is understood as a time-varying feedback, i.e., $u_t^{\star}(x)$ is, in general, non-homogeneous in $x$. This means that at any given time, the optimal corrective velocity experienced by one token embedding is different from another.  

The following Theorem \ref{thm:mb_transformer_sb} shows that the minimizing pair $t \mapsto (u_t^{\star},\rho_t^{\star})$ can be recovered by solving a system of coupled nonlinear integro-PDEs with bilinear boundary conditions. The proof for this result is given in the Supplementary Material Section \ref{sec:proof}.
\begin{theorem}[Schr\"{o}dinger System for Mean-field Token-steering]
\label{thm:mb_transformer_sb}
The minimizing tuple $(u_t^\star,\rho_t^\star)$ for problem \eqref{problem:spherical_control} is given by
\begin{align}
u_t^\star(x) &= \sigma \nabla \log \varphi_t(x), \label{eq:mb_transformer_optimal_control}\\
\rho_t^\star(x) &= \varphi_t(x)\widehat{\varphi}_t(x),   
\label{eq:schrodinger_factorization}
\end{align}
where $x \in \mathbb S^{d-1}$, and the positive real-valued functions $\varphi_t, \widehat \varphi_t$ on $\mathbb S^{d-1}$ satisfy 
\begin{subequations}
\label{eq:mb_transformer_schrodinger_system}
\begin{align}
\partial_t\varphi_t(x)  +  \big \langle  \nabla \varphi_t(x) &, P_x^\perp  (b_t[\rho_t^\star](x)) \big\rangle  +  \frac{\sigma}{2}  \Delta \varphi_t(x) \nonumber \\
&= - q_t[\rho_t^\star](x)\varphi_t(x) ,
\label{eq:mb_transformer_backward}
\\[1mm]
\partial_t\widehat{\varphi}_t(x)  +  \nabla\cdot \big( \widehat{\varphi}_t (x)& P_x^\perp( b_t[\rho_t^\star](x) )\big) -  \frac{\sigma}{2} \Delta\widehat{\varphi}_t(x) \nonumber  \\
&=  q_t[\rho_t^\star](x)\widehat{\varphi}_t(x),
\label{eq:mb_transformer_forward}
\end{align}
where $\langle \cdot, \cdot \rangle$ denotes the inner product induced by the round metric on $\mathbb S^{d-1}$, and \eqref{eq:mb_transformer_backward}-\eqref{eq:mb_transformer_forward} are subject to the multiplicative endpoint conditions
\begin{align}
\varphi_0(x)\widehat{\varphi}_0(x) &=  \rho_{\rm{initial}}(x), \label{eq:mb_transformer_endpoint_conditions-1} \\
\varphi_T(x)\widehat{\varphi}_T(x)  &= \rho_{\rm{final}}(x),\label{eq:mb_transformer_endpoint_conditions-2}
\end{align}
\end{subequations}
where 
\begin{align*}
b_t[\rho_t^\star](x) &= \frac{1}{{\mathcal Z_t[\rho_t^\star](x)}} \int_{\mathbb S^{d-1}}  e^{\beta\langle \mQ_t x,\mK_t y\rangle_{\mathbb R^d}} \,\mV_t \, y\,\rho_t^\star(y) \, {\rm d}y,\\
q_t[\rho_t^\star](x) &=  \int_{\mathbb S^{d-1}}   \widehat{\varphi}_t(y)  \frac{e^{\beta\langle \mQ_t y,\mK_t x\rangle_{\mathbb R^d}}}{\mathcal Z_t[\rho_t^\star](y)} \\
 & \quad \times \big\langle P_y^\perp  \left( \mV_tx-  b_t[\rho_t^\star](y) \right), \nabla\varphi_t(y) \big\rangle \,{\rm d}y,\\
\mathcal Z_t[\rho_t^\star](x) &= \int_{\mathbb S^{d-1}}  e^{\beta\langle \mQ_t x,\mK_t y\rangle_{\mathbb R^d}}  \rho_t^\star(y) \, {\rm d}y.
\end{align*}
\end{theorem}

\begin{remark}
The pair $\left(\varphi_t,\widehat{\varphi}_t\right)_{t\in[0,T]}$ is referred to as Schr\"{o}dinger factors \citep{caluya2021wasserstein} since they enable a factorization of $\rho_t^{\star}$; see \eqref{eq:schrodinger_factorization}.
\end{remark}

The system described in \eqref{eq:mb_transformer_schrodinger_system} is a mean-field generalization of the standard Schr\"odinger systems for the classical bridge problems; see \citet[Section~5]{chetrite2021schrodinger}, \citet{leonard2014survey}. It consists of nonlinear integro-PDEs that are coupled via the nonlocal attention transport $b_t[\rho_t^\star]$ and the reaction rate $q_t[\rho_t^\star]$. This structure differs from the standard Schr\"odinger bridge in that \eqref{eq:mb_transformer_backward}-\eqref{eq:mb_transformer_forward} then becomes a system of \emph{linear} advection-diffusion PDEs in $\left(\varphi_t,\widehat{\varphi}_t\right)$, see e.g., \citet[Theorem 2]{caluya2021wasserstein}. 

In the standard (i.e., non-mean-field) Schr\"odinger bridge, the linearity of \eqref{eq:mb_transformer_backward}-\eqref{eq:mb_transformer_forward} can be leveraged in computation via dynamic Sinkhorn recursions \citep[Section 8]{chen2021stochastic} over the function pair $\left(\varphi_{T},\widehat{\varphi}_0\right)$. In that case, the dynamic Sinkhorn (or iterative proportional fitting \citep{de2021diffusion}) recursion is known to be contractive with respect to Hilbert's projective metric \citep{bushell1973hilbert,georgiou2015positive}.  


In contrast, the nonlocal attention structure in $b_t,q_t$ makes \eqref{eq:mb_transformer_schrodinger_system} nonlinear, and it is not obvious if and how a Sinkhorn-like algorithm can be designed to compute the minimizing pair \eqref{eq:mb_transformer_optimal_control}-\eqref{eq:schrodinger_factorization}. Motivated by this, we propose a novel algorithm to obtain self-consistent solutions to \eqref{eq:mb_transformer_schrodinger_system}. 


\begin{algorithm}[!tp]
\caption{Mean-Field Token Steering Algorithm}
\label{alg:mfsb}
\begin{algorithmic}[1]
\Require $\beta,\mQ_t, \mK_t, \mV_t, \sigma, T,   \rho_{\rm initial},\rho_{\rm final}$, $N_K$, $N_{J}$, \texttt{tol}, $\theta \in(0,1]$.
\State Initialize  $\varphi_\cdot^{[0]},\widehat{\varphi}_\cdot^{[0]}$ and $\rho_\cdot^{(0)}=\varphi_\cdot^{[0]}\widehat{\varphi}_\cdot^{[0]}$.
\For{$k=0,\ldots N_K$} \State Compute $\mathcal Z_\cdot[\rho_\cdot^{(k)}],b_\cdot[\rho_\cdot^{(k)}]$.
\For{$j=0,\ldots,N_J$}
    \State Compute 
    \begin{align*}
        q_\cdot^{[j]}[\rho^{(k)}](x) &=\int_{\mathbb S^{d-1}}   \widehat{\varphi}_t^{[j]}(y)  \frac{e^{\beta\langle \mQ_t y,\mK_t x\rangle_{\mathbb R^d}}}{\mathcal Z_t[\rho_t^{(k)}](y)} \\
        &\times \Big\langle P_y^\perp  \left( \mV_tx-  b_t[\rho_t^{(k)}](y) \right), \nabla\varphi_t^{[j]}(y) \Big\rangle \,{\rm d}y. 
    \end{align*} 
        \State Integrate 
        \begin{align}
            \partial_t &\varphi_t^{[j+1]} (x)  + \big\langle  \nabla \varphi_t^{[j+1]} (x) , P_x^\perp  (b_t[\rho_t^{(k)}] (x)) \big\rangle  \nonumber \\
            &+  \frac{\sigma}{2}  \Delta \varphi_t^{[j+1]}(x) = - q_t^{[j]} [\rho^{(k)}] (x)\varphi_t^{[j+1]}(x) \label{eq:alg-back}
        \end{align}
        backward in time from $\varphi_T^{[j]}(x)$.
        \State Set
        \begin{align}
        \widehat \varphi_0^{[j+1]}(\cdot) = \frac{\rho_{\rm initial}(\cdot)}{\varphi_0^{[j+1]}(\cdot)}. \label{eq:alg-initial-impose}
        \end{align}
        \State Integrate 
        \begin{align}
           \partial_t &\widehat{\varphi}_t^{[j+1]} (x)  +  \nabla\cdot \big( \widehat{\varphi}^{[j+1]} (x)  P_x^\perp( b_t[\rho_t^{(k)}](x) )\big) \nonumber \\
           &-  \frac{\sigma}{2} \Delta\widehat{\varphi}_t^{[j+1]} (x) =  q_t^{[j]} [\rho^{(k)}] (x) \widehat{\varphi}_t^{[j+1]}(x) \label{eq:alg-for}
        \end{align}
        forward in time from $\widehat{\varphi}_0^{[j+1]}(x)$.
        \State Set
         \begin{align}
             \varphi_T^{[j+1]}(\cdot) = \frac{\rho_{\rm final}(\cdot)}{\widehat \varphi_T^{[j+1]}(\cdot)}. \label{eq:alg-final-impose}
         \end{align}
        \If{$
            \max_{t \in [0,T]} \max\Big\{ d_H\!\left( \varphi_t^{[j+1]}, \varphi_t^{[j]} \right)$,
            $ d_H \!\left( \widehat\varphi_t^{[j+1]}, \widehat\varphi_t^{[j]} \right) \Big\} < \texttt{tol}$}
            \State \textbf{break}
        \EndIf
    \EndFor
  \State Update
    \begin{align}
    \rho_\cdot^{(k+1)} = \theta\varphi_\cdot^{[j+1]}\widehat\varphi_\cdot^{[j+1]} +(1-\theta) \rho_\cdot^{(k)}.    \label{eq:iteration-nl2}
    \end{align}
    \If{$\max_{t \in [0,T]}  d_H\!\left(\rho_t^{(k+1)},\rho_t^{(k)}\right) <\texttt{tol} $}
    \State \textbf{break}
    \EndIf
    \State Set $\varphi_\cdot^{[0]} \gets \varphi_\cdot^{[j+1]}$ and $\widehat{\varphi}_\cdot^{[0]} \gets \widehat{\varphi}_\cdot^{[j+1]}$.
    \EndFor
    \State \Return $\varphi_\cdot^{[j+1]}$, $\widehat\varphi_\cdot^{[j+1]}$, and $\rho_\cdot^{(k+1)}$.
\end{algorithmic}
\end{algorithm}

\subsection{Computational Approach}
Building on Theorem \ref{thm:mb_transformer_sb}, we propose a generalized proportional fitting algorithm to compute the solution to the Schr\"odinger system \eqref{eq:mb_transformer_schrodinger_system}, enabling mean-field token steering.
Our computational methodology is detailed in Algorithm~\ref{alg:mfsb}, which comprises \emph{nested fixed-point iterations} indexed by $j$ and $k$. Specifically, the inner iteration index $j\in\{0,1,,\hdots,N_{J}\}$, and the outer iteration index $k\in\{0,1,,\hdots,N_{K}\}$. The initial guesses $\rho_\cdot^{(0)}$, $\varphi_\cdot^{[0]}$, and $\widehat{\varphi}_\cdot^{[0]}$ may be taken as the solution to the corresponding non-interacting Schr\"odinger bridge problem, obtained by setting the $b_t,q_t$ in \eqref{eq:mb_transformer_schrodinger_system} to zero.

In essence, the nested structure in Algorithm~\ref{alg:mfsb} handles the two sources of nonlinear coupling from the mean-field velocity $b_t$ and the reaction rate $q_t$ separately. Within the outer iteration that is indexed by $k$, the quantity $\rho_\cdot^{(k)}$ denotes the $k$-th approximation of the marginal flow. This approximation corresponds to mean-field velocity $b_\cdot[\rho_\cdot^{(k)}]$ and normalization $\mathcal Z_\cdot[\rho_\cdot^{(k)}]$. 
The inner iteration that is indexed by $j$ entails solving \eqref{eq:mb_transformer_backward}--\eqref{eq:mb_transformer_endpoint_conditions-2}, but with the fixed $b_t[\rho_t^{(k)}], \mathcal Z_t[\rho_t^{(k)}]$ and the reaction rate determined from the approximated marginal flow and delayed evaluation of the preceding Schr\"odinger factors.

\paragraph{Inner iteration} With these quantities fixed, the backward equation for $\varphi_t^{[j+1]}$ in \eqref{eq:alg-back} and the forward equation for $\widehat{\varphi}_t^{[j+1]}$ in \eqref{eq:alg-for} are \emph{linear} Kolmogorov equations. We use these linear equations for time marching, and update the marched terminal and initial values by imposing the prescribed final and initial marginal densities $\rho_{\rm final}, \rho_{\rm initial}$, as in \eqref{eq:alg-final-impose} and \eqref{eq:alg-initial-impose}, respectively.

\paragraph{Outer iteration} After the inner iteration completes, the outer loop updates $t \mapsto \rho_t^{(k)}$ according to \eqref{eq:iteration-nl2}, using a damped evaluation if found necessary for numerical stability. 
The nested iterations are repeated until the marginal flow and the Schr\"odinger pair are jointly self-consistent. 

\paragraph{Convergence criteria} We take our convergence criteria to be with respect to the Hilbert's projective metric defined for nonnegative functions $f,g$ in the Banach space $L^{\infty}\left(\mathbb{S}^{d-1}\right)$ as
\begin{align}
d_{H}(f,g):=\log\underset{x \in \Omega}{\sup}\dfrac{f(x)}{g(x)} - \log\underset{x\in \Omega}{\inf}\dfrac{f(x)}{g(x)}.
\label{defHilbertProjectiveMetric}    
\end{align}
The evaluation is restricted to the common support $\Omega \subseteq \mathbb S^{d-1}$ of the functions being compared.



\section{Numerical Experiments}\label{sec:Numerical}
We illustrate the proposed mean-field token steering framework on two experiments: one for $\mathbb{S}^{1}$ (circle) and the other for $\mathbb{S}^{2}$ (2-sphere in $\mathbb{R}^3$). In both experiments, we fix the noise strength parameter $\sigma=0.5$, inverse temperature $\beta=5$, and Algorithm \ref{alg:mfsb} parameters: ${\texttt{tol}}=10^{-6}$, $\theta=0.5$. In both, we perform a unimodal-to-bimodal steering over $[0,T]=[0,1]$ with
$$\rho_{\text{initial}}(x)
\propto
e^{\kappa_{0}\langle\mu_{0},x\rangle},
\quad
\rho_{\text{final}}(x)
\propto
\sum_{i=1}^{2}w_ie^{\kappa_{T}\langle\mu_{i},x\rangle},$$
for $\kappa_0,\kappa_T>0$, a convex weight vector $(w_1,w_2)$, and $\mu_0,\mu_1,\mu_2\in\mathbb{R}^{d}$ for $\mathbb{S}^{d-1}$. In other words, we steer a von Mises-Fisher distribution toward a two-component mixture of von Mises-Fisher distributions. We detail the other parameters and results for the specific experiments below.

\begin{figure*}[t!]
\centering
\includegraphics[width=\linewidth]{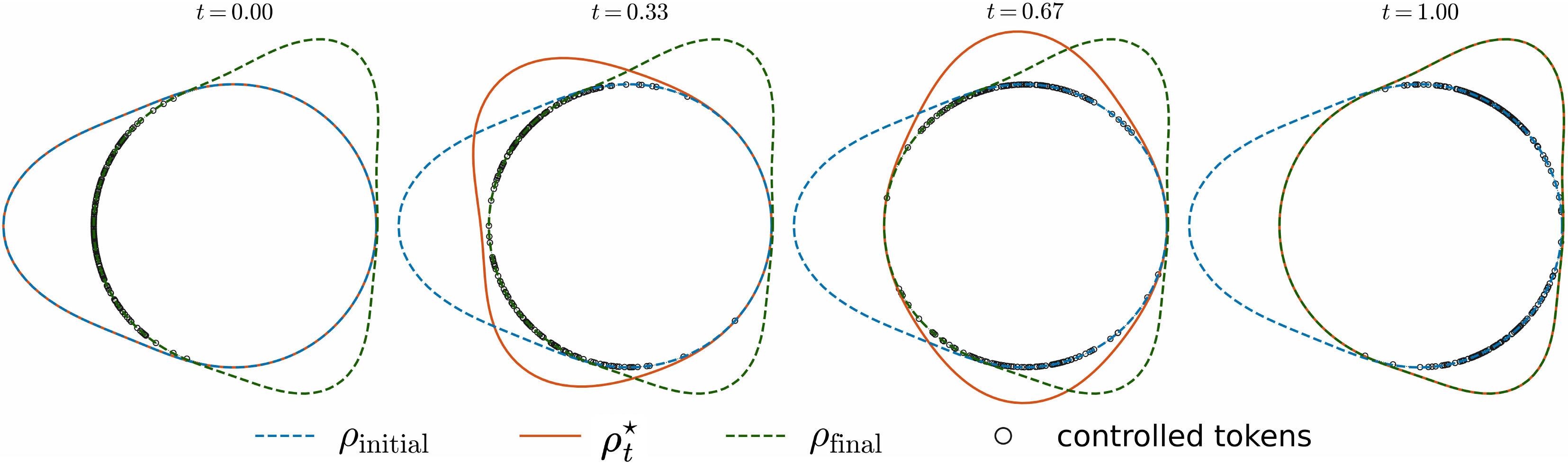}
\caption{Evolution of optimally controlled token distribution and tokens for the $\mathbb{S}^1$ experiment in Section \ref{subsec:CircleExperiment}.}
\label{fig:s1_snapshots}
\end{figure*}

\subsection{$\mathbb{S}^{1}$ Experiment}\label{subsec:CircleExperiment}

\paragraph{Setup} We use the following parameters for the initial and final distributions: $\kappa_0 = 6$, $\kappa_{T}=10$, $\mu_0 = (\cos\pi, \sin\pi)^{\top}$, $\mu_1 = (\cos(-1), \sin (-1))^{\top}$, $\mu_2=(\cos(1), \sin(1))^{\top}$, $w_1=0.4$, $w_2=0.6$. We use the self-attention parameters: $\bm{Q}_t=\bm{K}_t\equiv0.5\bm{I}_2$ (here $d=2$), and randomly generated
$$\bm{V}_t\equiv\begin{bmatrix}
    0.2180 & 0.2748\\
    0.0130 & 0.2177
\end{bmatrix}\quad\forall t\in[0,T].$$
In using Algorithm \ref{alg:mfsb} for this $\mathbb{S}^{1}$ experiment, we time-march the solutions for \eqref{eq:alg-back} and \eqref{eq:alg-for} using the Lie-Trotter splitting and $600$ uniform discretizations of the time horizon $[0,1]$. To approximate the spatial derivatives on the circle, we used the finite difference with $150$ discretizations of $[0,2\pi)$. To approximate the integrals in nonlocal terms, we used the periodic quadrature method. 

In applying Algorithm \ref{alg:mfsb} to this $\mathbb{S}^{1}$ experiment, we use $N_{J}=40$ inner iterations, and $N_{K} = 30$ outer iterations.  

\paragraph{Results} For the $\mathbb{S}^{1}$ numerical experiment setup described above, Figure~\ref{fig:s1_snapshots} shows four snapshots of the optimally controlled density $\rho_t^{\star}$ (solid curve) computed using Algorithm \ref{alg:mfsb}, together with $\rho_{\text{initial}},\rho_{\text{final}}$ (dashed curves). This plot demonstrates that Algorithm \ref{alg:mfsb} effectively steers the mean-field token distribution over the given horizon $[0,1]$. On these snapshots, we also superimposed the optimally controlled evolution of $200$ token embeddings (hollow circular markers), showing a higher concentration near the modes of $\rho_t^{\star}$, as expected. These token evolutions were computed via Euler-Maruyama simulation of the controlled dynamics \eqref{eq:dynamics-sphere-meanfield-controlled} with $u_t^{\star}$ obtained from Algorithm \ref{alg:mfsb} and with the same time discretization as before.
\begin{figure}[t!]
\centering
\includegraphics[width=\linewidth]{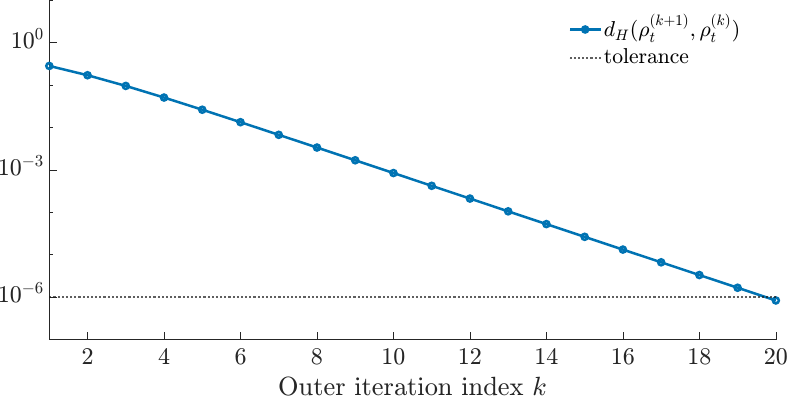}
\caption{Convergence of Algorithm \ref{alg:mfsb} for the experiment in Section \ref{subsec:CircleExperiment}.}
    \label{fig:s1_Hilbert}
\end{figure}

Figure~\ref{fig:s1_Hilbert} shows convergence of Algorithm \ref{alg:mfsb} with respect to the Hilbert metric \eqref{defHilbertProjectiveMetric}, for this $\mathbb{S}^1$ experiment.

An additional $\mathbb{S}^{1}$ experiment is reported in the Supplementary Material Section \ref{sec:AnotherS1Experiment}.


\subsection{$\mathbb{S}^{2}$ Experiment}\label{subsec:SphereExperiment}

\paragraph{Setup} Here, we use the following parameters for the initial and final distributions: $\kappa_0 = 30$, $\kappa_{T}=50$, $\mu_0 = (\sin(1)\cos(0.4), \sin(1)\sin(0.4),\cos(1))^{\top}$,
$$\mu_1 =
\frac{\mu_c\cos\alpha+\nu\sin\alpha}
{\|\mu_c\cos\alpha+\nu\sin\alpha\|_2}, \, \mu_2 = \frac{\mu_c\cos\alpha-\nu\sin\alpha}
{\|\mu_c\cos\alpha-\nu\sin\alpha\|_2},$$
$\mu_{c}:=-\mu_0$, $\nu$ is the unit tangent vector at $\mu_c$, $\alpha=0.9$ rad, $w_1=0.2$, $w_2=0.8$. For the self-attention parameters, we use $\beta=5$, fixed $\bm{Q}_t=\bm{K}_t\equiv\bm{I}_3$ (here $d=3$), and randomly generated
$$\bm{V}_t\equiv\begin{bmatrix}
0.2754 & 0.2554 & 0.0628 \\
0.3541 & 0.4465 & 0.1036 \\
0.1455 & 0.4481 & 0.0257
\end{bmatrix}\quad\forall t\in[0,T].$$
To apply Algorithm \ref{alg:mfsb} for this $\mathbb{S}^{2}$ experiment, we time-march the solutions for \eqref{eq:alg-back} and \eqref{eq:alg-for} using the Lie-Trotter splitting as before, with $200$ uniform discretizations of the time horizon $[0,1]$. To approximate the spatial derivatives on $\mathbb{S}^{2}$, we used the finite difference with $32\times 64$ discretizations of $[0,\pi]\times [0,2\pi)$. To approximate the integrals in nonlocal terms, we used the weighted quadrature method. 

In applying Algorithm \ref{alg:mfsb} to this $\mathbb{S}^{2}$ experiment, we use $N_{J}=50$ inner iterations, and $N_{K} = 60$ outer iterations.  

\begin{figure*}[t!]
\centering
\includegraphics[width=1\textwidth]{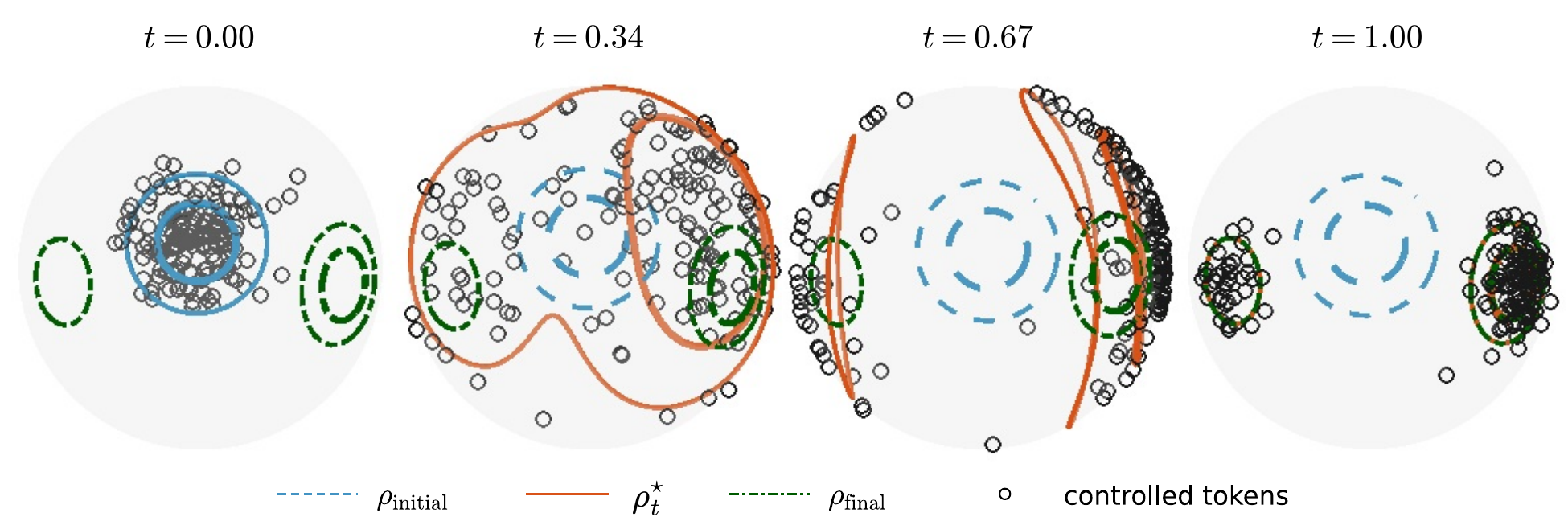}
\caption{Evolution of optimally controlled token distribution and tokens for the $\mathbb{S}^2$ experiment in Section \ref{subsec:SphereExperiment}.}
\label{fig:s2_snapshots}
\end{figure*}

\paragraph{Results} For the $\mathbb{S}^{2}$ numerical experiment setup described above, Figure~\ref{fig:s2_snapshots} shows four snapshots of the optimally controlled density $\rho_t^{\star}$ (solid contours) computed using Algorithm \ref{alg:mfsb}, together with $\rho_{\text{initial}},\rho_{\text{final}}$ (dashed contours), demonstrating the effectiveness of Algorithm \ref{alg:mfsb}. To avoid clutter, we suppressed the contour labels. The hollow circular markers in Figure~\ref{fig:s2_snapshots} show $160$ optimally controlled token embeddings, computed via Euler-Maruyama simulation of \eqref{eq:dynamics-sphere-meanfield-controlled} with $u_t^{\star}$ obtained from Algorithm \ref{alg:mfsb} and with the same time discretization as before.
\begin{figure}[t!]
\centering
\includegraphics[width=\linewidth]{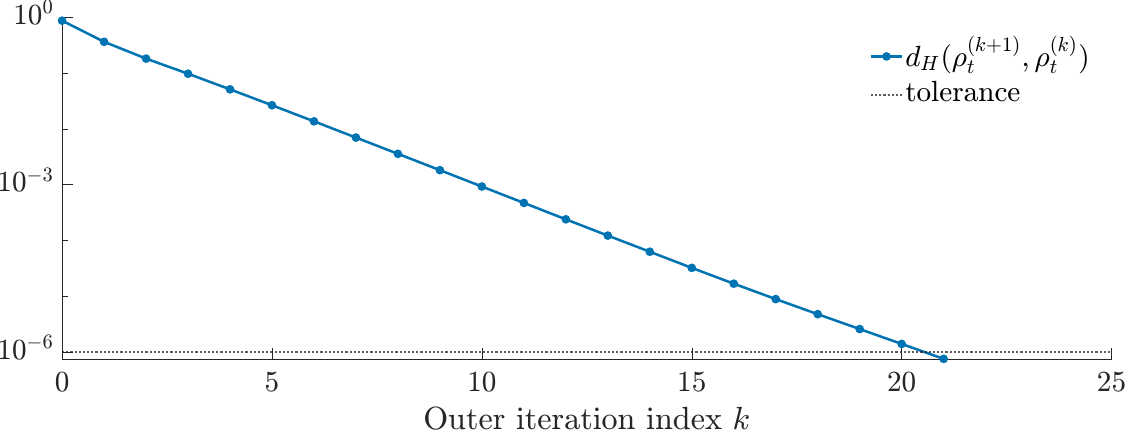}
\caption{Convergence of Algorithm \ref{alg:mfsb} for the experiment in Section \ref{subsec:SphereExperiment}.}
    \label{fig:s2_Hilbert}
\end{figure}

For this $\mathbb{S}^2$ experiment, Figure~\ref{fig:s2_Hilbert} shows convergence of Algorithm \ref{alg:mfsb} with respect to the Hilbert metric \eqref{defHilbertProjectiveMetric}.


\section{Concluding Remarks}\label{sec:Conclusions}
In this work, we have formulated the problem of steering self-attention-driven token dynamics in the mean-field limit with minimum control effort. The formulation, based on the theory of Schr\"odinger bridges, identifies the minimum-effort control as a feedback law that transports the token distribution between prescribed initial and final configurations over a finite transformer depth horizon. Our analysis characterizes the optimal controller via a forward--backward system of nonlinearly coupled integro-PDEs. To numerically solve such a system, we have proposed a generalized iterative proportional fitting algorithm adapted to the transformer mean-field setting.

Our primary contributions are theoretical: proposing the idea of mean-field token steering via control, connecting the formulation with the mean-field Schr\"{o}dinger bridge in the presence of attention-induced nonlocal drift, designing a variant of dynamic Sinkhorn recursion to solve the same, and demonstrating the feasibility of computation via simple numerical examples. While we empirically demonstrate local convergence, proving guarantees for the same remains a topic of ongoing research.

This work points towards multiple research avenues and connections within the rapidly evolving theory of transformer dynamics. 
Connecting the present result to recent analyses of representation collapse and meta-stable multi-clustering \citep{geshkovski2023emergence,bruno2026multiscale} is the most direct foreseeable research direction, and may provide a means to preserve the representation structure. Our numerical results support that possibility. 

Another connection is the recent line of research on activation steering (see, e.g., \citet{NEURIPS2023_81b83900,ICLR2025_e2dd5360,pmlr-v306-skifstad26a,pmlr-v306-adila26a}) that aims to guide an LLM during prompt context inference toward a desired semantic direction by modifying the model's internal representation at selected layers. The work here can be viewed as a distributional analogue of activation steering, where steering is achieved at the level of token-embedding populations. Therefore, an interesting direction for future work is to study how the controller obtained herein can be best realized (or approximated) by sparse interventions at finitely many layers, and how this restriction affects steering performance.

More broadly, extensions of the present framework to finite token populations, discrete layer depths, and multi-head attention indeed constitute an important step toward bridging the current theory with practical transformer architectures.


\subsection*{AI use statement}
During the preparation of this work, the authors used generative AI (OpenAI, Anthropic) to assist in developing the code for the numerical solvers and in post-processing of the figures. The authors reviewed, validated, and edited all the outputs and take full responsibility for the final contents, figures, claims, and artifacts presented here.

\balance
\bibliography{references}  

\clearpage
\onecolumn 

\setcounter{section}{0}
\setcounter{figure}{0}
\setcounter{table}{0}
\setcounter{equation}{0}
\renewcommand{\thesection}{\Alph{section}}
\renewcommand{\thefigure}{S\arabic{figure}}
\renewcommand{\thetable}{S\arabic{table}}
\renewcommand{\theequation}{S\arabic{equation}}

\begin{center}
    \vbox{
        \hsize\textwidth
        \linewidth\textwidth
        \hrule height 4pt depth 0pt width \textwidth
        \vspace{4mm}
        \centering
        {\Large\bf Supplementary Material: Mean-field Token Steering \par}
        \vspace{4mm}
        \hrule height 1pt depth 0pt width \textwidth
        \vspace{5mm}
    }
\end{center}


\section{Proof of Theorem \ref{thm:mb_transformer_sb}}\label{sec:proof}
This section provides the proof of Theorem \ref{thm:mb_transformer_sb}.
We start by writing the Lagrangian of Problem~\eqref{problem:spherical_control} as
\begin{align} 
\mathcal L(u_\cdot,\rho_\cdot,\Lambda_\cdot) &= \frac{1}{2\sigma}\int_0^T \int_{\mathbb S^{d-1}} \|u_t(x)\|^2 \rho_t(x)\, {\rm d} x \, {\rm d} t \nonumber \\
&+\int_0^T \int_{\mathbb S^{d-1}} \Lambda_t(x) \Bigg[ \partial_t\rho_t(x)   +\nabla \cdot \Big( \rho_t(x) \Big(P_x^\perp \Big( b_t[\rho_t](x) \Big) + u_t(x) \Big) \Big) - \frac{\sigma }{2}\Delta \rho_t(x) \Bigg]  \, {\rm d} x \, {\rm d} t, \label{eq:mean_field_augmented_lagrangian} 
\end{align}
where $\Lambda_t: \mathbb S^{d-1} \to \mathbb R$ is a Lagrange multiplier that enforces the dynamical constraint \eqref{eq:spherical_FP}.

Using integration by parts, we write
\begin{align} 
\int_0^T \int_{\mathbb S^{d-1}} \Lambda_t (x)\,\partial_t\rho_t (x) \, {\rm d} x \, {\rm d} t 
&= \int_{\mathbb S^{d-1}} \Big\{\Lambda_T(x)\rho_T(x) -\Lambda_0(x)\rho_0(x)\Big\} \,{\rm d} x  \nonumber \\
&- \int_0^T \int_{\mathbb S^{d-1}} \rho_t(x) \partial_t \Lambda_t(x) \, {\rm d} x \, {\rm d} t. \label{eq:mean_field_time_IBP} 
\end{align}
Next, define 
\begin{equation} \label{eq:v-t}
v_t(x) := P_x^\perp \Big( b_t[\rho_t](x) \Big) +  u_t(x), \qquad \forall x \in \mathbb S^{d-1}.
\end{equation} 
By the divergence theorem \citep[Theorem~16.48]{lee2012introduction}, we have 
\begin{align*} 
\int_{\mathbb S^{d-1}} &\Lambda_t(x) \nabla \cdot \Big( \rho_t(x) v_t(x) \Big) {\rm d} x = - \int_{\mathbb S^{d-1}} \rho_t \big\langle v_t(x), \nabla\Lambda_t(x) \big\rangle {\rm d} x,
\end{align*}
since $\mathbb S^{d-1}$ is a closed manifold, i.e., $\partial\mathbb S^{d-1}=\varnothing$. 
Hence, substituting \eqref{eq:v-t} into the previous equation yields 
\begin{align}
\int_{\mathbb S^{d-1}} \Lambda_t(x) \nabla  \cdot \Big(\rho_t(x) \Big(P_x^\perp \Big( b_t[\rho_t](x) \Big)  +  u_t(x) \Big) \Big) {\rm d} x = - \int_{\mathbb S^{d-1}} \rho_t(x) \left\langle P_x^\perp \Big( b_t[\rho_t](x)\Big) + u_t(x) , \nabla\Lambda_t(x) \right\rangle {\rm d} x, \label{eq:mean_field_drift_IBP_lagrangian} 
\end{align}
recalling that $\langle \cdot, \cdot \rangle$ denotes the inner product induced by the round metric on $\mathbb S^{d-1}$.

By Green's second identity \citep[Page~51]{lee2018introduction}, we also obtain that
\begin{equation}
-\frac{\sigma}{2} \int_{\mathbb S^{d-1}} \Lambda_t (x) \, \Delta \rho_t(x) \,{\rm d} x = -\frac{\sigma}{2}\int_{\mathbb S^{d-1}} \rho_t (x) \, \Delta \Lambda_t(x) \,{\rm d} x. \label{eq:mean_field_diffusion_IBP_lagrangians} 
\end{equation} 
By \eqref{eq:mean_field_time_IBP}, \eqref{eq:mean_field_drift_IBP_lagrangian}, and \eqref{eq:mean_field_diffusion_IBP_lagrangians}, the Lagrangian in \eqref{eq:mean_field_augmented_lagrangian} can be rewritten as
\begin{align} 
&\mathcal L(u_\cdot,\rho_\cdot,\Lambda_\cdot) = \int_{\mathbb S^{d-1}} \Big\{ \Lambda_T(x)\rho_T(x)  -\Lambda_0(x)\rho_0(x) \Big\}  \,{\rm d} x \nonumber \\
&+ \int_0^T \int_{\mathbb S^{d-1}} \rho_t(x) \Bigg[ \frac{\|u_t(x)\|^2}{2\sigma}- \partial_t \Lambda_t(x)  - \left \langle P_x^\perp \Big( b_t[\rho_t](x) \Big) +  u_t(x), \nabla\Lambda_t(x) \right\rangle  - \frac{\sigma}{2} \Delta \Lambda_t(x)  \Bigg]  \, {\rm d} x \, {\rm d} t.
\label{eq:mean_field_reduced_lagrangian} 
\end{align}

\subsection*{Optimization over the controller $t\mapsto u_t$}
The first-order optimality condition for $u_\cdot^\star$ is such that the G\^{a}teaux derivative 
\begin{align}
  \frac{\rm d}{\rm d \eta}  \mathcal L (u_\cdot^\star + \eta  \, \delta u_\cdot^\star, \rho_\cdot^\star, \Lambda_\cdot) \Big\vert_{\eta=0} =0,
\end{align}
for every perturbation $t \mapsto \delta u_t^\star$, with $\delta u_t^\star(x) \in T_x \mathbb S^{d-1}$, $x \in \mathbb S^{d-1}$. This implies 
\begin{align*}
\frac1\sigma\int_0^T \int_{\mathbb S^{d-1}} \rho_t^\star(x) \Big\langle u_t^\star(x) - \sigma \nabla\Lambda_t(x), \delta u_t^\star(x)
\Big\rangle \,{\rm d} x\,{\rm d} t =0.
\end{align*}
That is, for almost every $t \in [0,T]$,
\begin{align}
    u_t^\star = \sigma \nabla \Lambda_t , \qquad \rho_t^\star\text{-almost everywhere} . \label{eq:optimal_control_mean_field}
\end{align}

\subsection*{Optimization over the marginal flow $t \mapsto \rho_t$}
Let us substitute \eqref{eq:optimal_control_mean_field} into \eqref{eq:mean_field_reduced_lagrangian}. This yields the following evaluation at the optimizer $t \mapsto (\rho_t^\star, u_t^\star)$. 
\begin{align}
    \mathcal L(\rho_\cdot^\star&,\Lambda_\cdot) = \int_{\mathbb S^{d-1}} \Big\{ \Lambda_T(x)\rho_T^\star(x)  -\Lambda_0(x)\rho_0^\star(x) \Big\}  \,{\rm d} x\nonumber \\
    &+ J(\rho_\cdot^\star,\Lambda_\cdot)+ \int_0^T \int_{\mathbb S^{d-1}} \Bigg\{ \rho_t^\star(x) \Bigg[- \partial_t \Lambda_t(x) - \frac{\sigma}{2} \Vert \nabla\Lambda_t (x) \Vert^2 - \frac{\sigma}{2} \Delta \Lambda_t(x) \Bigg] \Bigg\} \, {\rm d} x \, {\rm d} t, \label{eq:Lagrnagian-only-rho-star}
\end{align}
where
\begin{align}\label{eq:J-star}
    &J(\rho_\cdot^\star, \Lambda_\cdot) := - \int_0^T \int_{\mathbb S^{d-1}} \rho_t^\star(x) \Big\langle P_x^\perp \Big(b_t[\rho_t^\star](x) \Big), \nabla\Lambda_t(x) \Big\rangle \, {\rm d} x\,{\rm d} t.
\end{align}

We let $ t \mapsto \delta \rho_t^\star$ be a perturbation of the optimal PDF-valued curve $t \mapsto \rho_t^\star$, so that
\begin{align*}
    \int_{\mathbb S^{d-1}} \delta \rho_t^\star =0.
\end{align*}
Further, since the endpoint optimal PDFs $\rho_0^\star,\rho_T^\star$ are fixed as in \eqref{eq:spherical_endpoints}, 
\begin{align}\label{eq:rho-star-0-T}
    \delta \rho_0^\star = \delta\rho_T^\star=0.
\end{align}
Then, the first-order optimality condition with respect to $\rho_\cdot$ is such that
\begin{align}
  \frac{\rm d}{\rm d \eta}  \mathcal L (\rho_\cdot^\star + \eta  \, \delta \rho_\cdot^\star, \Lambda_\cdot) \Big\vert_{\eta=0} =0.
\end{align}
Therefore, by \eqref{eq:Lagrnagian-only-rho-star}--\eqref{eq:rho-star-0-T}, we have
\begin{align}\label{eq:optimality-before-subs}
- \int_0^T \int_{\mathbb S^{d-1}} \delta \rho_t^\star(x) \Bigg[\partial_t \Lambda_t(x) +  \frac{\sigma}{2} \Vert \nabla \Lambda_t (x) \Vert^2 + \frac{\sigma}{2} \Delta \Lambda_t(x) \Bigg] \, {\rm d}x \, {\rm d}t +   \frac{\rm d}{\rm d\eta}J(\rho_\cdot^\star+\eta  \, \delta \rho_\cdot^\star, \Lambda_\cdot)\Big\vert_{\eta=0} = 0. 
\end{align}

Indeed,
\begin{align}
\frac{\rm d}{\rm d\eta}J(\rho_\cdot^\star+\eta \, \delta \rho_\cdot^\star, \Lambda_\cdot)  &= - \int_0^T \int_{\mathbb S^{d-1}} \Bigg\{ \delta \rho_t^\star(x) \left\langle P_x^\perp  \Big(b_t[\rho_t^\star+\eta \delta \rho_t^\star](x)\Big) , \nabla\Lambda_t(x) \right\rangle  \nonumber \\ 
&\qquad+ \Big (\rho_t^\star(x) + \eta \, \delta \rho_t^\star (x) \Big) \frac{\rm d}{\rm d\eta}  \Big \langle P_x^\perp \Big(b_t[\rho_t^\star+\eta \delta \rho_t^\star](x) \Big),
\nabla\Lambda_t(x) \Big\rangle  \Bigg\}  \, {\rm d} x \, {\rm d} t. \label{eq:drift_variation_product_rule}
\end{align}
Notice that
\begin{align}
\frac{\rm d}{\rm d \eta}  \left\langle P_x^\perp \Big( b_t[\rho_t^\star+\eta  \, \delta \rho_t^\star](x) \Big),\nabla\Lambda_t(x)\right\rangle &=\left\langle P_x^\perp \Big(
\frac{\rm d}{\rm d \eta} b_t[\rho_t^\star+ \eta  \, \delta \rho_t^\star](x) \Big), \nabla\Lambda_t(x) \right\rangle. \label{eq:eta-intermediate}
\end{align}
Therefore, \eqref{eq:drift_variation_product_rule} and \eqref{eq:eta-intermediate} lead  to
\begin{align}
\frac{\rm d}{\rm d \eta}J(\rho_\cdot^\star+\eta \, \delta \rho_\cdot^\star, \Lambda_\cdot)\Big|_{\eta=0} & = -  \int_0^T \int_{\mathbb S^{d-1}} \Bigg\{ \delta \rho_t^\star(x) \left\langle P_x^\perp \Big( b_t[\rho_t^\star+\eta \, \delta \rho_t^\star]\Big|_{\eta=0}(x) \Big),  \nabla\Lambda_t(x) \right\rangle \nonumber\\
&\quad+ \rho_t^\star(x) \left\langle P_x^\perp \Big(\frac{\rm d}{\rm d \eta} b_t[\rho_t^\star + \eta \, \delta \rho_t^\star]\Big|_{\eta=0}(x)\Big), \nabla\Lambda_t(x) \right\rangle \Bigg \} \, {\rm d} x \, {\rm d} t. 
\label{eq:drift_variation_before_eta_zero}
\end{align}

Recall that
\begin{align}
    b_t[\rho_t^\star](x) &= \frac{1}{{\mathcal Z_t[\rho_t^\star](x)}} \int_{\mathbb S^{d-1}}  e^{\beta\langle \mQ_t x,\mK_t y\rangle_{\mathbb R^d}} \,\mV_t \, y\,\rho_t^\star(y) \, {\rm d}y,\\
    \mathcal Z_t[\rho_t^\star](x) &= \int_{\mathbb S^{d-1}}  e^{\beta\langle \mQ_t x,\mK_t y\rangle_{\mathbb R^d}}  \rho_t^\star(y) \, {\rm d}y. \label{eq:z-t}
\end{align}
Let us set
\begin{align} \label{eq:y-t}
    \mathcal Y_t[\rho_t^\star] (x) := \int_{\mathbb S^{d-1}}  e^{\beta\langle \mQ_t x,\mK_t y\rangle_{\mathbb R^d}} \,\mV_t \, y\,\rho_t^\star(y) \, {\rm d}y.
\end{align}
Hence,
\begin{align*}
     b_t[\rho_t^\star] &= \frac{\mathcal Y_t[\rho_t^\star]}{\mathcal Z_t[\rho_t^\star]}, \nonumber\\
      b_t[\rho_t^\star+ \eta \, \delta \rho_t^\star] &= \frac{\mathcal Y_t[\rho_t^\star+ \eta \, \delta \rho_t^\star]}{\mathcal Z_t[\rho_t^\star+ \eta \, \delta \rho_t^\star]}.
\end{align*}

Note that by using \eqref{eq:z-t} and \eqref{eq:y-t}, we can rewrite the previous equation as
\begin{align*}
    b_t[\rho_t^\star+ \eta \, \delta \rho_t^\star] &= \frac{\mathcal Y_t[\rho_t^\star]+ \eta \, \mathcal Y_t[\delta \rho_t^\star]}{\mathcal Z_t[\rho_t^\star] +\eta \,  \mathcal Z_t[\delta \rho_t^\star]}.
\end{align*}
Thus,
\begin{align}
    \left.  \frac{\rm d}{\rm d \eta} b_t[\rho_t^\star+\eta \delta \rho_t^\star]\right|_{\eta=0}  
& = \frac{1}{\mathcal Z_t[\rho_t^\star]} \Big(\mathcal Y_t[\delta \rho_t^\star] - b_t[\rho_t^\star]\mathcal Z_t[\delta \rho_t^\star]\Big), \nonumber\\
   \left.  \frac{\rm d}{\rm d \eta} b_t[\rho_t^\star+\eta \delta \rho_t^\star]\right|_{\eta=0} (x) &=  \frac{1}{{\mathcal Z_t[\rho_t^\star](x)}} \int_{\mathbb S^{d-1}}  e^{\beta\langle \mQ_t x,\mK_t y\rangle_{\mathbb R^d}} \Big( \,\mV_t \, y - b_t[\rho_t^\star](x) \Big) \, \delta\rho_t^\star(y) \, {\rm d}y. \nonumber
\end{align}
Additionally,
\begin{align}
   P_x^\perp \Big(   \left.  \frac{\rm d}{\rm d \eta} b_t[\rho_t^\star+\eta \delta \rho_t^\star]\right|_{\eta=0} (x) \Big)&=  \frac{1}{{\mathcal Z_t[\rho_t^\star](x)}} \int_{\mathbb S^{d-1}}  e^{\beta\langle \mQ_t x,\mK_t y\rangle_{\mathbb R^d}}  P_x^\perp  \Big( \,\mV_t \, y - b_t[\rho_t^\star](x) \Big) \, \delta\rho_t^\star(y) \, {\rm d}y, \label{eq:d-eta-d-b-t}
\end{align}
where we have used the fact that $\mathcal Z_t[\rho_t^\star](x)$ and $e^{\beta\langle \mQ_t x,\mK_t y\rangle_{\mathbb R^d}}$ and scalars for any fixed $x$.

Using \eqref{eq:d-eta-d-b-t}, we can rewrite \eqref{eq:drift_variation_before_eta_zero} as
\begin{align}
\frac{\rm d}{\rm d \eta}J&(\rho_\cdot^\star+\eta \, \delta \rho_\cdot^\star, \Lambda_\cdot)\Big|_{\eta=0}  = -  \int_0^T \int_{\mathbb S^{d-1}} \delta \rho_t^\star(x) \left\langle P_x^\perp \Big(  b_t[\rho_t^\star](x) \Big),  \nabla\Lambda_t(x) \right\rangle   \, {\rm d} x \, {\rm d} t. \nonumber\\
& -\int_0^T \int_{\mathbb S^{d-1}}  \frac{\rho_t^\star(x)}{\mathcal Z_t[\rho_t^\star](x)} \left \langle \int_{\mathbb S^{d-1}}  e^{\beta\langle \mQ_t x,\mK_t y\rangle_{\mathbb R^d}}  P_x^\perp  \Big( \,\mV_t \, y - b_t[\rho_t^\star](x) \Big) \, \delta\rho_t^\star(y) \, {\rm d}y, \nabla\Lambda_t(x) \right\rangle \, {\rm d} x \, {\rm d} t. \label{eq:d-J-d-eta-1}
\end{align}
By exchanging the variables $x$ and $y$ in the second integral on the right-hand side of \eqref{eq:d-J-d-eta-1}, we obtain
\begin{align}
\frac{\rm d}{\rm d \eta}J(\rho_\cdot^\star+\eta \, \delta \rho_\cdot^\star&, \Lambda_\cdot)\Big|_{\eta=0}  = -  \int_0^T \int_{\mathbb S^{d-1}}  \Bigg\{ \delta \rho_t^\star(x) \Bigg[\left\langle P_x^\perp \Big(  b_t[\rho_t^\star](x) \Big),  \nabla\Lambda_t(x) \right\rangle \nonumber\\
& + \int_{\mathbb S^{d-1}}\frac{\rho_t^\star(y) e^{\beta\langle \mQ_t y,\mK_t x\rangle_{\mathbb R^d}}}{\mathcal Z_t[\rho_t^\star](y)} \left \langle    P_y^\perp  \Big( \,\mV_t \, x - b_t[\rho_t^\star](y) \Big)  , \nabla\Lambda_t(y) \right\rangle \, {\rm d} y \Bigg] \Bigg\} \, {\rm d}x\, {\rm d} t. \label{eq:d-J-d-eta-2}
\end{align}

By substituting \eqref{eq:d-J-d-eta-2} into \eqref{eq:optimality-before-subs}, we obtain 
\begin{align*}
   \int_0^T \int_{\mathbb S^{d-1}} \Bigg\{ \delta \rho_t^\star(x)  \Big[ \partial_t \Lambda_t(x) &+  \frac{\sigma}{2} \Vert \nabla\Lambda_t (x) \Vert^2 + \frac{\sigma}{2} \Delta \Lambda_t(x) \\
   &+ \left\langle P_x^\perp \Big(  b_t[\rho_t^\star](x) \Big),  \nabla\Lambda_t(x) \right\rangle +q_t[\rho_t^\star] (x) \Big] \Bigg\} \, {\rm d} x \, {\rm d}t =0.
\end{align*}
Since $\int_{\mathbb S^{d-1}} \delta \rho_t^\star(x) =0$, we deduce that the quantity in brackets must be spatially independent. That is,
\begin{align}
     \partial_t \Lambda_t(x) &+  \frac{\sigma}{2} \Vert \nabla \Lambda_t (x) \Vert^2 + \frac{\sigma}{2} \Delta \Lambda_t(x) + \left\langle P_x^\perp \Big(  b_t[\rho_t^\star](x) \Big),  \nabla\Lambda_t(x) \right\rangle +q_t[\rho_t^\star] (x) = f(t)
\end{align}
where
\begin{align}
q_t[\rho_t^\star](x) :=  \int_{\mathbb S^{d-1}}  \frac{ \rho_t^\star(y)e^{\beta\langle \mQ_t y,\mK_t x\rangle_{\mathbb R^d}}}{\mathcal Z_t [\rho_t^\star](y)} \left\langle P_y^\perp
\Big(\mV_t \, x-b_t[\rho_t^\star](y)\Big),\nabla\Lambda_t(y)\right\rangle\,{\rm d} y, \label{eq:mean_field_correction_y}
\end{align}
for some $f: [0,T] \to \mathbb R^d$.

Since $f$ depends only on time, we can absorb it into $\Lambda_t$ via the transformation $$\Lambda_t + \int_0^t f(s) \, {\rm d} s \to \Lambda_t$$ and write the first-order optimality in $\rho_\cdot ^\star$ as the Hamilton--Jacobi--Bellman equation 
\begin{align}
    \partial_t \Lambda_t(x) +  \frac{\sigma}{2} \Vert \nabla\Lambda_t (x) \Vert^2 + \frac{\sigma}{2} \Delta \Lambda_t(x) + \left\langle P_x^\perp \Big(  b_t[\rho_t^\star](x) \Big),  \nabla\Lambda_t(x) \right\rangle = -q_t[\rho_t^\star] (x).
\end{align}

Therefore,  we can summarize the optimality conditions for Problem~\eqref{problem:spherical_control} as follows.
\begin{subequations}  \label{eq:mean_field_optimality_hopf_cole}
\begin{align} 
u_t^\star(x) &= \sigma\, \nabla\Lambda_t(x), \label{eq:control-lambda}\\
\partial_t\rho_t^\star(x)  &+\nabla \cdot \left(\rho_t^\star\Big(P_x^\perp \Big( b_t[\rho_t^\star](x)\Big)+ u_t^\star(x) \Big)\right)= \frac\sigma2\Delta\rho_t^\star(x),\\
 \partial_t \Lambda_t(x) &+  \frac{\sigma}{2} \Vert \nabla \Lambda_t (x) \Vert^2 + \frac{\sigma}{2} \Delta \Lambda_t(x) + \left\langle P_x^\perp \Big(  b_t[\rho_t^\star](x) \Big),  \nabla\Lambda_t(x) \right\rangle = -q_t[\rho_t^\star] (x), \label{eq:mean_field_optimality_HJB}\\
\rho_{0}^\star(x) &= \rho_{\rm initial} (x), \qquad \rho_{T}^\star(x) = \rho_{\rm final}(x),  \label{eq:BC-Lambda}
\end{align}
\end{subequations}
where $q_t[\rho^\star_t]$ is as in \eqref{eq:mean_field_correction_y}.

\subsection*{The forward-backward dynamics}
Next, we introduce a change of variable
\begin{equation}
\Lambda_t(x)=\log \varphi_t(x),
\label{eq:hopf-cole}
\end{equation}
where $\varphi_t$ is a positive function for all $t \in [0,T]$. This transformation is known as the \emph{Cole--Hopf} \citep{cole1951quasi,hopf1950partial,teter2025hopf} or \emph{Fleming logarithmic transformation} \citep{fleming2005logarithmic}. By \eqref{eq:hopf-cole}, Eq.~\eqref{eq:control-lambda} can be written as
\begin{align*}
    u_t^\star(x) = \sigma \, \nabla \log \varphi_t(x),
\end{align*}
recovering \eqref{eq:mb_transformer_optimal_control}.

By \eqref{eq:hopf-cole}, \eqref{eq:mean_field_optimality_HJB}, and the identity
\begin{equation*}
\left\|  \nabla\Lambda_t(x)  \right\|_{g}^2 + \Delta\Lambda_t(x)= \frac{ \Delta\varphi_t(x) }{\varphi_t(x)},
\end{equation*}
we can obtain the dynamics of $\varphi_t$ as in \eqref{eq:mb_transformer_backward}.

Introducing a second nonnegative function $\widehat{\varphi}_t$ through the factorization
\begin{equation*}
    \rho_t^\star(x) =  \varphi_t(x)\widehat{\varphi}_t(x),
\end{equation*}
as in \eqref{eq:schrodinger_factorization}, 
we deduce that $\widehat{\varphi}_t$ satisfies \eqref{eq:mb_transformer_forward}.
Lastly, by \eqref{eq:schrodinger_factorization}, we can obtain \eqref{eq:mb_transformer_endpoint_conditions-1} and \eqref{eq:mb_transformer_endpoint_conditions-2}.

Combining the preceding optimality conditions and the defined transformations yields the Schr\"odinger system stated in the theorem. 

\section{Additional $\mathbb{S}^{1}$ Experiment}\label{sec:AnotherS1Experiment}
We report an additional mean-field token-steering experiment over $\mathbb{S}^{1}$ that supplements our results in Section \ref{subsec:CircleExperiment}. Here, we focus on comparing the uncontrolled versus optimally controlled evolution of the token distribution in the mean-field (i.e., number of tokens $N\rightarrow\infty$) limit.

\begin{figure*}[t!]
\centering
\includegraphics[width=\linewidth]{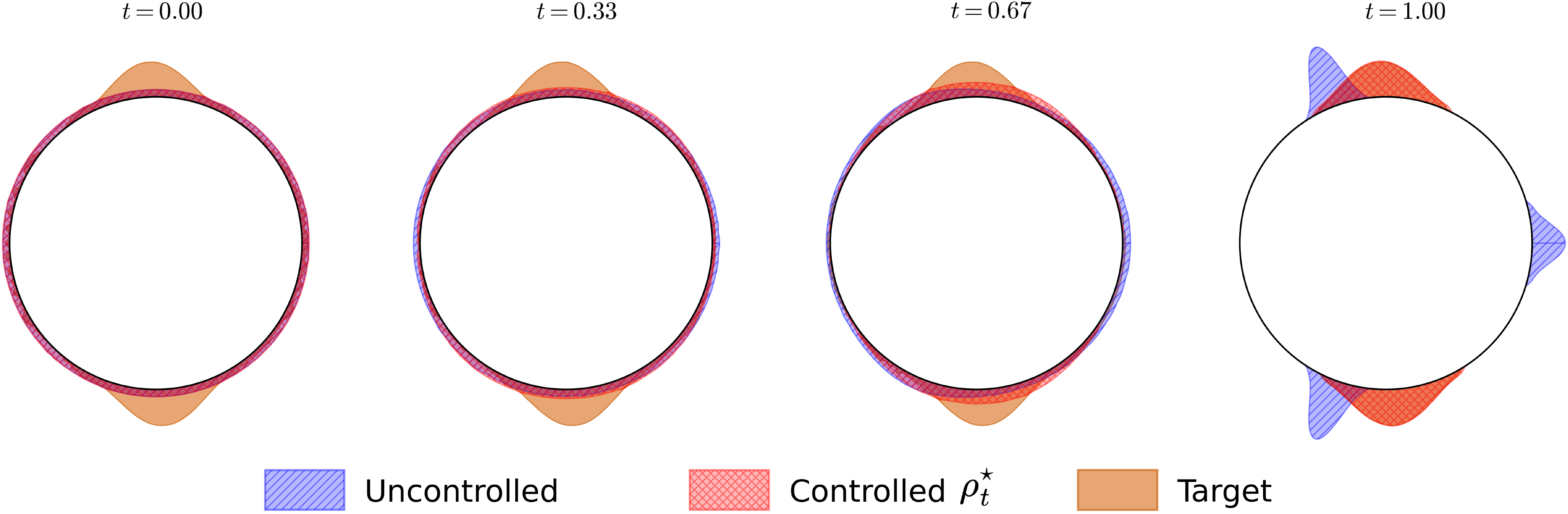}
\caption{Evolution of uncontrolled versus optimally controlled token distribution for the $\mathbb{S}^1$ experiment in Supplementary Section \ref{sec:AnotherS1Experiment}.}
\label{fig:s1_ControlledVersusUncontrolled}
\end{figure*}

Inspired by the experiment reported in \citet[Fig. 4]{Bruno2025b}, we set $\rho_{\text{initial}}$ to be uniform on $\mathbb{S}^1$ slightly perturbed by white noise. It is known that for this $\rho_{\text{initial}}$, the uncontrolled mean-field dynamics for certain choice of parameters gives rise to a $3$ mode meta-stable clustering. To investigate the possibility of countering this trend via control, we set $\rho_{\text{final}}$ as an equal-weight mixture of two von Mises-Fisher densities on $\mathbb{S}^1$, given by
$$\rho_{\text{final}}\propto \sum_{i=1}^{2}w_ie^{\kappa_{T}\langle\mu_{i},x\rangle},$$
with $w_1=w_2=0.5$, $\kappa_{T}=25$, $\mu_i=(\cos\theta_i,\sin\theta_i)^{\top}$ for $i\in\{1,2\}$, $\theta_1:=1.6072$ rad, and $\theta_2 :=\theta_1 + \pi$. In Fig. \ref{fig:s1_ControlledVersusUncontrolled}, this $\rho_{\text{final}}$ is shown as ``target". We use the parameters: $\sigma=0.2$, $T=1$, $\beta=5$, $\bm{Q}_t=\bm{K}_t\equiv\bm{I}_2$ (here $d=2$), and randomly generated
$$\bm{V}_t\equiv\begin{bmatrix}
    0.0305 & 0.1123\\
    0.1171 & 0.0885
\end{bmatrix}\quad\forall t\in[0,T].$$
We follow the space-time discretization as in \citet[Section 5]{Bruno2025b}, and in Fig. \ref{fig:s1_ControlledVersusUncontrolled}, compare the uncontrolled (i.e., with prior mean-field token dynamics) versus controlled (i.e., with corrected mean-field token dynamics) evolution of the token densities over $t\in[0,1]$. The controlled evolution is computed using Algorithm \ref{alg:mfsb} proposed here with the same algorithm parameters as in Section \ref{subsec:CircleExperiment}. Fig. \ref{fig:s1_ControlledVersusUncontrolled} shows that the optimal control accomplishes the uniform to target bimodal token steering, overcoming the natural tendency of forming $3$ modes over the same horizon.

\end{document}